\documentclass[a4paper,10pt]{amsart}
\usepackage{amsmath,amssymb,amsxtra,amsfonts,amsthm,url,comment,microtype,enumerate}

\usepackage{a4wide}

\theoremstyle{plain}
\newtheorem{theorem}{Theorem}
\newtheorem{lemma}{Lemma}
\newtheorem{corollary}{Corollary}

\begin{document}

\title{The Calkin--Wilf problems on $\ell$-ary partition functions}

\subjclass[2020]{Primary 05A17; Secondary 05A15, 05A16, 11P82}
\keywords{Partitions, Regular Sequences, Mahler Functions, Probability Measures}

\author{Michael Coons}

\address{Aarhus University, Ny Munkegade 118, 8000 Aarhus C, Denmark, and \newline \hspace*{.3cm} California State University, 400 West First Street, Chico, California 95929, USA}
\email{mjc@math.au.dk, mjcoons@csuchico.edu}

\date{\today}

\begin{abstract} In their 1998 precursor to their celebrated paper ``Recounting the rationals,'' Calkin and Wilf listed nine open problems related to the function $b_{k,\ell}(n)$, the number of partitions of $n$ into powers of $\ell$ with each part appearing at most $k$ times. In this short paper, we answer five of their eight remaining open problems. As a consequence of our results, we extend and generalise results of both Carlitz and Reznick, and answer an open question of Reznick.
\end{abstract}

\maketitle

\section{Introduction}

Let $k,\ell\geqslant 2$ and let $b_{k,\ell}(n)$ denote the number of partitions of $n$ into powers of $\ell$ with each power appearing at most $k$ times. That is, $b_{k,\ell}(n)$ equals the number of solutions to $$n=a_0+a_1\ell+a_2\ell^2+\cdots+a_m\ell^n+\cdots,$$ with $a_m\in\{0,1,\ldots,k\}$ for each $m\geqslant 0$. The study of $\ell$-ary partitions is old and has an extensive literature---to much to recount in such a short paper. Some highlights of the asymptotic theory include the 1940 result of Mahler \cite{M1940} who showed that $\log b_{\infty,2}(n)\sim (\log n)^2/\log 4$ and the more complete understanding of $b_{\infty,\ell}(n)$ by de Bruijn \cite{dB1948} in 1948. We refer the interested reader to the binary partition focussed paper of Reznick \cite{R1990} for extensive references and details on $b_{k,2}(n)$, and for more general partition theory to the volumes of Andrews \cite{A1998} and Andrews and Eriksson \cite{AE2004}. Our study starts with the work of Calkin and Wilf. In their paper, ``Recounting the rationals''~\cite{CW2000}, they showed that the sequence $(b_{2,2}(n)/b_{2,2}(n+1))_{n\geqslant 0}$ of ratios of binary partition functions enumerates the nonnegative rational numbers, already in reduced form and without repeats. This $(3+\varepsilon)$-page gem started as a nineteen-page preprint \cite{CWpre}, released in 1998 and updated in 2009, that contained several more results and musings as well as nine open problems for further research---we answered their ninth question in 2014 in joint work with Tyler \cite{CT2014}. More formally, after noting that {\em it is an easy exercise to show that} $b_{\ell^2-1,\ell}(n)=\lfloor n/\ell\rfloor +1$, Calkin and Wilf state the following eight (remaining) open problems.

\begin{quote}
\begin{enumerate}
\item[CW1.] Is $k=\ell^2-1$ the first time we get smooth growth in $b_{k,\ell}(n)$? That is, if $\ell\leqslant k<\ell^2-1$, then the growth is not smooth?
\item[CW2.] Is $k=\ell^2-1$ the first time we get linear maximum growth?
\item[CW3.] Is $k=\ell^2-1$ the first time we get linear minimum growth?
\item[CW4.] Is $k=\ell^2-1$ the first time we get non-constant monotone growth?
\item[CW5.] Is the sequence $b_{k,\ell}(n)$ ever log-convex?
\item[CW6.] What can be said about consecutive digit sets?
\item[CW7.] What can be said of general digit sets?\item[CW8.] For which pairs $(k,\ell)$ is $b_{k,\ell}(n)$ monotonic in $n$?
\end{enumerate}
\end{quote}

\noindent In this paper, we answer five of the eight remaining questions:~CW1, CW3, CW4, CW5 and CW8. Our proofs use some old, and some not-so-old, asymptotic results on the generating functions of these partitions. In particular, we exploit the $\ell$-regularity of the sequences $b_{k,\ell}(n)$ to great benefit.

\section{The Calkin--Wilf problems on $\ell$-ary partition functions}

For finite $k$, one can immediately write down the generating function of $b_{k,\ell}(n)$ as the infinite product $$B_{k,\ell}(z):=\sum_{n\geqslant 0}b_{k,\ell}(n)z^n=\prod_{j\geqslant 0}\big(1+z^{\ell^j}+z^{2\ell^j}+\cdots+z^{k\ell^j}\big).$$ Note that \begin{equation}\label{eq:Bmahler}B_{k,\ell}(z)-(1+z+z^2+\cdots+z^k)B_{k,\ell}(z^\ell)=0,\end{equation} which implies that $B_{k,\ell}(z)$ is an $\ell$-Mahler function. Thus, a result of ours with Bell~\cite[Thm.~1]{BC2017} on the asymptotics of Mahler functions applies to give \begin{equation}\label{eq:B}B_{k,\ell}(z)=C_{k,\ell}(z)(1-z)^{-\log_\ell(k+1)}(1+o(1)),\end{equation} as $z\to 1^-$, where $\log_\ell$ denotes the base-$\ell$ logarithm, $C_{k,\ell}(z)$ is a positive (possibly) oscillatory term, which in $(0,1)$ is bounded away from $0$ and $\infty$, is real analytic and satisfies $C_{k,\ell}(z)=C_{k,\ell}(z^\ell)$. 

With this asymptotic in hand, we are ready for our first result, which answers CW3 affirmatively. 

\begin{theorem} The function $b_{k,\ell}(n)$ grows at least linearly if and only if $k\geqslant \ell^2-1$.
\end{theorem}

\begin{proof} It is clear that $b_{k,\ell}(n)$ is a monotone increasing function of $k$. Since we know that $b_{\ell^2-1,\ell}(n)$ grows linearly, to prove our result, it is enough to show that $b_{k,\ell}(n)$ is not bounded below by a linear function of $n$ for any $k<\ell^2-1$. To this end, suppose that there is a real number $\alpha>0$ such that $\alpha\cdot n\leqslant b_{\ell^2-2,\ell}(n)$. Then, as $z\to1^-$, for some $\beta>0$, we have $$\frac{\alpha z}{(1-z)^2}=\sum_{n\geqslant 0}\alpha nz^n\leqslant B_{k,\ell}(z)\leqslant\frac{\beta}{(1-z)^{\log_\ell(k+1)}},$$ where the right-most inequality follows from our above-mentioned result with Bell. But this implies that $2<\log_\ell(k+1)$, so that $k\geqslant \ell^2-1$, which proves the result.
\end{proof}

Considering the sentiment of question CW3, we offer the following result.

\begin{theorem}\label{thm:k=ellm} If $m\geqslant 1$ and $k\in(\ell^m-1,\ell^{m+1}-1)$, then there are constants $\alpha_1,\alpha_2>0$ such that for $n$ large enough $\alpha_1 n^{m-1}\leqslant b_{k,\ell}(n)\leqslant \alpha_2 n^m$.
\end{theorem}

\begin{proof} As well, appealing to the fact that $b_{k,\ell}(n)$ is a monotone increasing function of $k$, the result will follow provided $b_{\ell^m-1,\ell}(n)\sim c_mn^{m-1}$ for some constant $c_m>0$. To see this, we start by using the polynomial identity \begin{equation}\label{eq:polyID}1+z+z^2+\cdots+z^k=\frac{1-z^{k+1}}{1-z},\end{equation} to obtain, for any $m\geqslant 1$, that $$B_{\ell^m-1,\ell}(z)=\prod_{j\geqslant 0}\left(\frac{1-z^{\ell^{m+j}}}{1-z^{\ell^j}}\right)=\prod_{j=0}^{m-1}\big(1-z^{\ell^j}\big)^{-1}=(1-z)^{-m}\prod_{j=0}^{m-1}\big(1+z+z^{2}+\cdots+z^{\ell^j-1}\big)^{-1}.$$ Since $z=1$ is a pole of order $m$ of $B_{\ell^m-1,\ell}(z)$ and $$\lim_{z\to 1}\prod_{j=0}^{m-1}\big(1+z+z^{2}+\cdots+z^{\ell^j-1}\big)^{-1}=\ell^{m(m-1)/2},$$ a standard result in analytic combinatorics---see Flajolet and Sedgewick \cite[p.~256, IV.26]{FSbook}---gives that \begin{equation}\label{eq:c_m}b_{\ell^m-1,\ell}(n)=\frac{n^{m-1}}{\ell^{m(m-1)/2}(m-1)!}\left(1+O\left(\frac{1}{n}\right)\right),\end{equation} which proves the result.
\end{proof}

The proof of Theorem \ref{thm:k=ellm}---in particular, Equation \eqref{eq:c_m}---gives the following corollary regarding a special case of CW1.

\begin{corollary}\label{cor:mono} For each $m\geqslant 1$, the function $b_{\ell^m-1,\ell}(n)$ is smooth.
\end{corollary}

Before returning to smoothness and CW1, we address CW4 and CW8 concerned with monotonicity and CW5 on log-convexity. In particular, we prove the following theorem, which provides the complete answer to CW8, as well as a negative answer to CW4.

\begin{theorem}\label{thm:mono} The function $b_{k,\ell}(n)$ is monotone if and only if $\ell$ divides $k+1$. Moreover, $b_{k,\ell}(n)$ is never strictly monotone.
\end{theorem}

\begin{proof} Suppose that $\ell$ does not divide $k+1$. The functional equation \eqref{eq:Bmahler} implies \begin{align*}\sum_{n\geqslant 0} b_{k,\ell}(n)z^n&=\sum_{n\geqslant 0}b_{k,\ell}(n)z^{\ell n}\big(1+z+z^2+\cdots+z^k\big) \\
&=\sum_{n\geqslant 0}b_{k,\ell}(n)\big(z^{\ell n}+z^{\ell n+1}+\cdots+z^{\ell n+k}\big)\\
&=\sum_{n=0}^{k-1}b_{k,\ell}(n)z^n+\sum_{n\geqslant k}\left(\sum_{j=\lceil {(n-k)}/{\ell}\rceil}^{\lfloor n/{\ell}\rfloor}b_{k,\ell}(j)\right)z^n.\end{align*} Using this, we compare the values of $b_{k,\ell}(k+j\ell)$ and $b_{k,\ell}(k+1+j\ell)$ for each $j\geqslant 0$. We have $$b_{k,\ell}(k+j\ell)=\sum_{m=j}^{\lfloor k/\ell\rfloor+j}b_{k,\ell}(m)>\sum_{m=j+1}^{\lfloor k/\ell\rfloor+j}b_{k,\ell}(m)=b_{k,\ell}(k+1+j\ell),$$ where we have used the fact that $\lfloor (k+1)/\ell\rfloor=\lfloor k/\ell\rfloor$ since $\ell$ does not divide $k+1$.

Now, suppose that $\ell$ divides $k+1$, and write $k=a\ell-1$. Using the polynomial identity $$1+z+z^2+\cdots+z^k=(1+z+\cdots+z^{\ell-1})(1+z^{\ell}+z^{2\ell}+\cdots+z^{a\ell}),$$ we have $$B_{k,\ell}(z)=\prod_{j\geqslant 0}\big(1+z^{\ell^j}+z^{2\ell^j}+\cdots+z^{(\ell-1)\ell^j}\big)\cdot \prod_{j\geqslant 1}\big(1+z^{\ell^j}+z^{2\ell^j}+\cdots+z^{a\ell^j}\big),$$ where the first product is equal to the series $\sum_{n\geqslant 0}z^n$ since its coefficients are the number of ways to write $n$ as the sum of powers of $\ell$ with each power used at most $\ell-1$ times---this is just the standard base-$\ell$ expansion, which is unique. The second product can be written as $\sum_{n\geqslant 0}g(n)z^n$, where, since the index $j\geqslant 1$, $g(n)=0$ if $n$ is not divisible by $\ell$, and $g(\ell n)>0$ for all $n\geqslant 0$. We thus have that the coefficients $b_{k,\ell}(n)$ of $B_{k,\ell}(z)$ satisfy \begin{equation}\label{eq:bg}b_{k,\ell}(n)=\sum_{m\leqslant n}g(n),\end{equation} which is monotonically increasing, since $g(n)\geqslant 0$, though not strictly, since $g(n)=0$ for infinitely many $n$.
\end{proof}

Our next result will rely heavily on the fact that $b_{k,\ell}(n)$ is an $\ell$-regular sequence. Indeed, this is a direct consequence of the functional equation \eqref{eq:Bmahler}; see Becker \cite[Thm.~2]{B1994}. Since $b_{k,\ell}(n)$ is an $\ell$-regular sequence, a classical result of Allouche and Shallit \cite[Thm.~2.2(e)]{AS1992} gives that there are an integer $d\geqslant 1$ (taken to be minimal), matrices ${\bf B}_0,{\bf B}_1,\ldots,{\bf B}_{\ell-1}\in\mathbb{Z}^{d\times d}$ and a vector ${\bf v}\in\mathbb{Z}^{d\times 1}$ such that $$b_{k,\ell}(n)={\bf e}_1^T {\bf B}_{i_0}{\bf B}_{i_1}\cdots {\bf B}_{i_s}{\bf v},$$ where ${\bf e}_1$ is the first $d\times 1$ standard basis vector and $(n)_\ell=i_s\cdots i_1i_0$ is the base-$\ell$ expansion of $n$. Set ${\bf B}:={\bf B}_0+{\bf B}_1+\cdots+{\bf B}_{\ell-1}$. Then, the asymptotics in \eqref{eq:B} imply that the spectral radius $\rho:=\rho({\bf B})=k+1$. The maximal growth of $b_{k,\ell}(n)$ can also be quantified. As we saw above, for each $m\geqslant 1$, there is a constant $c_m>0$ such that $b_{\ell^m-1,\ell}(n)\sim c_m n^{m-1}$. Since $b_{k,\ell}(n)$ is an increasing sequence in $k$, we have for $k\in(\ell^m-1,\ell^{m+1}-1)$ and $n$ large enough $$(c_m/2)n^{m-1}\leqslant b_{k,\ell}(n)\leqslant 2c_{m+1}n^m.$$ The maximal order of growth of an $\ell$-regular sequence is $\log_\ell\rho^*$, where $\rho^*:=\rho^*(\{{\bf B}_0,{\bf B}_1,\ldots,{\bf B}_{\ell-1}\})$ is the joint spectral radius of the associated finite set of matrices \cite{C2017}. Thus, since $k+1>\ell^m$, $$\log_\ell\rho^*\leqslant m<\log_\ell(k+1)=\log_\ell\rho.$$ In particular, $\rho^*<\rho$. It is convenient now to record a result of Dumas \cite{D2013}---since $\rho^*<\rho$, we have \begin{equation}\label{eq:dumas}\sum_{n\leqslant x}b_{k,\ell}(n)=\psi_{k,\ell}(\log_\ell x)x^{\log_\ell(k+1)}(1+o(1)),\end{equation} where $\psi_{k,\ell}(x)$ is a positive $\log_\ell(\rho/r)$-H\"older continuous $1$-periodic function for any $r\in(\rho^*,\rho)$. Combining this result with the last part of the proof of Theorem \ref{thm:mono}, we are able to generalise and extend a result of Reznick \cite[Thm.~4.7]{R1990}.

\begin{theorem} Suppose that $\ell$ divides $k+1$ and write $k+1=a\ell$ for some $a\in\mathbb{N}$. Then $$b_{k,\ell}(n)=\frac{\psi_{a,\ell}(\log_\ell n)}{a+1}\cdot n^{\log_\ell(a+1)}(1+o(1)).$$ where $\psi_{a,\ell}(x)$ is the positive H\"older continuous $1$-periodic function from \eqref{eq:dumas}.
\end{theorem}

\begin{proof} Let $g(n)$ be the function defined in the proof of Theorem \ref{thm:mono}, so that $$G(z):=\sum_{n\geqslant 0}g(n)z^n=\prod_{j\geqslant 1}\big(1+z^{\ell^j}+z^{2\ell^j}+\cdots+z^{a\ell^j}\big)=B_{a,\ell}(z^\ell).$$ Now, Dumas's result implies that $$\sum_{m\leqslant n}b_{a,\ell}(m)=\psi_{a,\ell}(\log_\ell n)n^{\log_\ell(a+1)}(1+o(1)).$$ But, the coefficients $g(m)$ of $B_{a,\ell}(z^\ell)$ satisfy $$g(m)=\begin{cases} 0, & \mbox{if $\ell$ does not divide $m$};\\ b_{a,\ell}(m/\ell),& \mbox{if $\ell$ does divide $m$}.\end{cases}$$ Using this with \eqref{eq:bg} gives \begin{align*}b_{k,\ell}(n)=\sum_{m\leqslant n}g(n)&=\sum_{m\leqslant n/\ell}b_{a,\ell}(m)\\
&=\psi_{a,\ell}(\log_\ell (n/\ell))(n/\ell)^{\log_\ell(a+1)}(1+o(1))\\
&=\frac{\psi_{a,\ell}(\log_\ell n)}{a+1}\cdot n^{\log_\ell(a+1)}(1+o(1)),\end{align*} where we have used that $\psi_{a,\ell}(\log_\ell (n/\ell))=\psi_{a,\ell}(\log_\ell n)$ since $\psi_{a,\ell}(x)$ is $1$-periodic.
\end{proof}

The proof of Theorem \ref{thm:mono} also allows us to immediately answer question CW5 concerning log-convexity. Recall that a sequence $a(n)$ is {\em log-convex} provided $$a(n)^2\leqslant a(n-1)\cdot a(n+1)$$ for all $n$. We say the sequence is {\em eventually log-convex} if the inequality holds for $n$ large enough.

\begin{theorem} For each $k\neq \ell-1$, the sequence $b_{k,\ell}(n)$ is not eventually log-convex.
\end{theorem}

\begin{proof} If $k<\ell-1$, then $b_{k,\ell}(n)$ takes the values $0$ and $1$ infinitely many times. For any one of the infinitely many times $b_{k,\ell}(n)=1$ and $b_{k,\ell}(n+1)=0$, we have that $$b_{k,\ell}(n)^2=1>0=b_{k,\ell}(n-1)\cdot b_{k,\ell}(n+1).$$ 

Now suppose that $k>\ell-1$. Then, Theorem \ref{thm:mono} implies that $b_{k,\ell}(n)$ is not strictly monotonic. Since $b_{k,\ell}(n)$ is not strictly monotonic, there are infinitely many triples $(n-1,n,n+1)$ such that $$b_{k,\ell}(n-1)<b_{k,\ell}(n)\quad\mbox{and}\quad b_{k,\ell}(n)\geqslant b_{k,\ell}(n+1),$$ and so for each of these $n$, $b_{k,\ell}(n)^2>b_{k,\ell}(n-1)\cdot b_{k,\ell}(n+1)$. 
\end{proof}

Having addressed monotonicity and log-convexity, we return to smoothness. We already have a partial result towards this, Corollary \ref{cor:mono}, which says that $b_{k,\ell}(n)$ has smooth growth for $k+1=\ell^m$ for any $m\geqslant 1$. It turns out, these are the only values of $k$ for which $b_{k,\ell}(n)$ grows smoothly.

\begin{theorem}\label{thm:smooth} The function $b_{k,\ell}(n)$ is smooth if and only if $k+1=\ell^m$ for some $m\geqslant 1$. 
\end{theorem}

To prove Theorem \ref{thm:smooth} we must better understand the oscillatory function $C_{k,\ell}(z)$ from \eqref{eq:B}. For a given Mahler function, this can be difficult, but, in our specific case, we are able to use (again, again) the polynomial identity \eqref{eq:polyID}. We combine this with an old result of de Bruijn to understand when these functions oscillate. De Bruijn \cite{dB1948} studied the asymptotics for general $\ell$ when $k=\infty$. In particular, an immediate consequence of his main result is, as $z\to1^-$, that \begin{equation}\label{eq:Tell} B_{\infty,\ell}(z)^{-1}=\prod_{j\geqslant 0}(1-z^{\ell^j})=C_\ell(z)(1-z)^{1/2}\ell^{-\log_\ell^2(1-z)/2}(1+o(1)),\end{equation} where $C_\ell(z)$ is a positive oscillatory term, which in $(0,1)$ is bounded away from $0$ and $\infty$, is real analytic and satisfies $C_\ell(z)=C_\ell(z^\ell)$. We stress that, here, $C_\ell(z)$ is truly oscillatory (so non-constant), which is why we have left off the `(possibly)' in this case; indeed, de Bruijn \cite[Eq.~(2.15)]{dB1948}, making the convenient substitution $z=e^{-s}$ showed that \begin{equation}\label{eq:Cell} C_\ell(e^{-s})=\exp\left(\sum_{n\in\mathbb{Z}} a(n)s^{-\frac{2\pi i n}{\log\ell}}\right)\end{equation} with \begin{equation}\label{eq:an}a(n)=\begin{cases}\big(\gamma_2-\gamma^2/2+\pi^2/12+(\log\ell)^2/12\big)/\log\ell&\mbox{if $n=0$},\\\Gamma\big({2\pi i n}/{\log\ell}\big)\zeta\big(1+{2\pi i n}/{\log\ell}\big)/\log \ell&\mbox{if $n\neq 0$},\end{cases}\end{equation} where $\gamma$ is the Euler--Mascheroni constant, $\gamma_2$ is the second Stieltjes constant, and $\Gamma(z)$ and $\zeta(z)$ are the standard Euler and Riemann functions, respectively. Relating this to $B_{k,\ell}(z)$, using \eqref{eq:polyID} and \eqref{eq:Tell}, as $z\to1^-$, we have \begin{equation}\label{eq:BTM}B_{k,\ell}(z)=\prod_{j\geqslant 0}\left(\frac{1-z^{(k+1)\ell^j}}{1-z^{\ell^j}}\right)=\frac{C_\ell(z^{k+1})}{C_\ell(z)}\cdot\frac{\sqrt{k+1}}{\ell^{\log_\ell^2(k+1)/2}}\cdot(1-z)^{-\log_\ell(k+1)}(1+o(1)).\end{equation} Comparing equations \eqref{eq:B} and \eqref{eq:BTM}, we have that \begin{equation}\label{eq:C=CC}C_{k,\ell}(z)=\frac{C_\ell(z^{k+1})}{C_\ell(z)}\cdot\frac{\sqrt{k+1}}{\ell^{\log_\ell^2(k+1)/2}}.\end{equation} 

To exploit this relationship, we use the following lemma.

\begin{lemma}\label{lem:Cklconst} If $b_{k,\ell}(n)$ is smooth, then $C_{k,\ell}(z)$ is constant.
\end{lemma}

\begin{proof} Suppose that $b_{k,\ell}(n)$ is smooth. Then, there is a $c_{k,\ell}>0$ such that \begin{equation}\label{eq:basymp}b_{k,\ell}(n)=c_{k,\ell}n^{\log_\ell(k+1)-1}+o\big(n^{\log_\ell(k+1)-1}\big).\end{equation} The exponent on $n$ comes from the fact that $\sum_{n\leqslant x}b_{k,\ell}(n)\asymp n^{\log_\ell(k+1)}$ and that, if $b_{k,\ell}(n)$ is smooth, its growth is determined by the average growth. Equation \eqref{eq:basymp} implies that as $z\to 1^-$, $$B_{k,\ell}(z)=c_{k,\ell}{\rm Li}_{1-\log_\ell(k+1)}(z)+o\big({\rm Li}_{1-\log_\ell(k+1)}(z)\big)=\frac{c_{k,\ell}\cdot \Gamma(\log_\ell(k+1))}{(1-z)^{\log_\ell(k+1)}}(1+o(1)),$$ where ${\rm Li}_s(z)$ is the standard polylogarithm. Equating this asymptotic with \eqref{eq:B}, as $z\to 1^-$, $${C_{k,\ell}(z)}(1+o(1))={c_{k,\ell}\cdot \Gamma(\log_\ell(k+1))}(1+o(1)),$$ and so, since $C_{k,\ell}(z)=C_{k,\ell}(z^\ell)$, we have that $C_{k,\ell}(z)=c_{k,\ell}\cdot \Gamma(\log_\ell(k+1))$ identically.
\end{proof}

Lemma \ref{lem:Cklconst} immediate implies, via \eqref{eq:C=CC}, that if $b_{k,\ell}(n)$ is smooth, then $C_{\ell}(z)=\gamma C_\ell(z^{k+1})$ for some constant $\gamma>0$. Indeed, we can say a bit more.

\begin{lemma}\label{lem:osc} If there is an integer $k+1\geqslant 2$ and a real $\gamma$ such that $C_\ell(z)=\gamma C_\ell(z^{k+1})$, then $\gamma=1$.
\end{lemma}

\begin{proof} First, since $C_\ell(z)$ is positive and bounded in $(0,1)$, $\sup_{z\in(0,1)}C_\ell(z)$ exists and is positive. Further, since $z\mapsto z^{k+1}$ is surjective on $(0,1)$, if $C_\ell(z)=\gamma C_\ell(z^{k+1})$, we have that $$\sup_{z\in(0,1)}C_\ell(z)=\sup_{z\in(0,1)}\gamma C_\ell(z^{k+1})=\gamma\sup_{z\in(0,1)}C_\ell(z),$$ so that $\gamma=1$. 
\end{proof}

We are now ready to prove Theorem \ref{thm:smooth}.

\begin{proof}[Proof of Theorem \ref{thm:smooth}] Suppose $b_{k,\ell}(n)$ is smooth. By Lemmas \ref{lem:Cklconst} and \ref{lem:osc}, we have $C_\ell(z)=C_{\ell}(z^{k+1})$. Consider the sum $$\psi(s):=\sum_{n\in\mathbb{Z}}a(n)s^{-\frac{2\pi i n}{\log\ell}}$$ from \eqref{eq:Cell}. Since $C_\ell(e^{-s})=\exp(\psi(s))$ and $C_\ell(z)=C_{\ell}(z^{k+1})$, we have \begin{align*}0=\psi(s)-\psi(s(k+1))&=\sum_{n\in\mathbb{Z}\backslash\{0\}} a(n)s^{-\frac{2\pi i n}{\log\ell}}\big(1-e^{-{2\pi i n\log_\ell(k+1) }})\\
&=\sum_{n\in\mathbb{Z}\backslash\{0\}} a(n)\big(1-e^{-{2\pi i n\log_\ell(k+1) }})e^{-\frac{2\pi i n \log s}{\log\ell}}.\end{align*} This equality holds for all positive $s$, so applying this equality for $s=\ell^{\vartheta}$ with $\vartheta\in[0,2\pi)$, gives $$\sum_{n\in\mathbb{Z}\backslash\{0\}} a(n)\big(1-e^{-{2\pi i \log_\ell(k+1) n}})e^{-{2\pi i n \vartheta}}=0$$ for all $\vartheta\in[0,2\pi)$. Thus, by Cantor's uniqueness theorem \cite{C1870} (see also Ash \cite{A1989}), for each $n\in\mathbb{Z}\backslash\{0\}$, $$a(n)\big(1-e^{-{2\pi i \log_\ell(k+1) n}})=0.$$ Since $\Gamma(it)\zeta(1+it)$ is nonzero for real $t$, we have $e^{-2\pi i \log_\ell(k+1)n}=1$ for each $n$, and so $\log_\ell(k+1)\in\mathbb{Z}$. That is, $k+1$ is a power of $\ell$, which completes the proof.
\end{proof}

Theorem \ref{thm:smooth} is the complete generalisation of a result Carlitz \cite{C1965}, who, encouraged by Bateman, proved that $b_{5,2}(n)$ is not smooth. Theorem \ref{thm:smooth} provides a generalised answer to a related question of Reznick \cite[p.~474]{R1990} on binary partition functions; that is, $$\liminf_{n\to\infty} b_{k,\ell}(n)n^{-\log_\ell(k+1)-1}=\limsup_{n\to\infty} b_{k,\ell}(n)n^{-\log_\ell(k+1)-1}$$ if and only if $k=\ell^m-1$ for some $m\geqslant 1$. 

\section{Further remarks}

There now remain three open Calkin--Wilf problems on $\ell$-ary partition functions. 

Questions CW6 and CW7 are pretty general---these ask about consecutive or general digit sets. What we do know is that the generating function, say $P(z)$, of either of these is still an infinite product of the form $P(z)=\prod_{j\geqslant 0}p\big(z^{\ell^j}\big)$ where $p(z)$ is a polynomial with coefficients in $\{0,1\}$. Thus, we still obtain the asymptotic $P(z)=C(z)(1-z)^{-\log_\ell p(1)}(1+o(1))$ as $z\to 1^-$, as in \eqref{eq:B}. To get more specific properties may be difficult---one doesn't immediately have polynomial identities like the ones we've heavily relied on. But, the coefficients are still regular sequences, so there is certainly some structure. We find it an interesting question whether or not it is possible to know if the coefficients of an infinite product such as $P(z)$ are smooth or monotone by only considering properties of the generating function. For a specified polynomial $p(z)$ one can probably do something, but are there general properties of $p(z)$ that determine smoothness or monotonicity?

We find question CW2 extremely interesting, and even more so, understanding the maximal growth of any regular sequence. In general, these questions have to do with being able to determine the joint spectral radius $\rho^*$ of the related set of matrices, a question which is well-known in many mathematical areas, as well as engineering, and still extremely open. For the specific question CW2, it would be enough to show that the maximal growth of $B_{\ell^2-2,\ell}(n)$ is strictly sublinear. As a generalisation, we ask whether it is always the case that $$\limsup_{n\to\infty}\frac{\log b_{k,\ell}(n)}{\log n}< \limsup_{n\to\infty}\frac{\log b_{k+1,\ell}(n)}{\log n},$$ for all $k\geqslant \ell-1$? What about the limit inferiors? 

\subsection*{Acknowledgements}

We thank Simon Kristensen for helpful conversations. This work was supported by Villum Fonden research grant VIL83378 and by the David W.~and Helen E.~F.~Lantis Endowment.

\bibliographystyle{amsplain}
\providecommand{\bysame}{\leavevmode\hbox to3em{\hrulefill}\thinspace}
\providecommand{\href}[2]{#2}


\end{document}